\documentclass[12pt]{article}%
\usepackage{amsmath,amssymb, amsthm, epsfig}
\usepackage{hyperref}
\usepackage{amsmath}
\usepackage{amssymb}
\usepackage{hyperref}
\usepackage{color}
\usepackage{ulem}
\usepackage{epsfig}
\usepackage{amsfonts}
\usepackage{graphicx}%
\providecommand{\U}[1]{\protect\rule{.1in}{.1in}}
\newlength{\hchng}
\newlength{\vchng}
\newtheorem{theorem}{Theorem}

\newtheorem{proposition}[theorem]{Proposition}
\newtheorem{corollary}[theorem]{Corollary}
\theoremstyle{definition}
\newtheorem{definition}[theorem]{Definition}

\theoremstyle{remark}

\begin{document}

\date{}
\author{Daniel Pellegrino $\&$ Eduardo V. Teixeira}
\title{Norm optimization problem for linear operators in classical Banach spaces}
\maketitle

\begin{abstract}
The main result of the paper shows that, for $1<p<\infty$ and
$1\leq q<\infty,$ a linear operator
$T\colon\ell_{p}\rightarrow\ell_{q}$ attains its norm if, and only
if, there exists a not weakly null maximizing sequence for $T$
(counterexamples can be easily constructed when $p=1$). For
$1<p\not =q<\infty,$ as a consequence of the previous result we
show that any not weakly null maximizing sequence for a norm
attaining operator $T:\ell_{p}\rightarrow\ell_{q}$ has a
norm-convergent subsequence (and this result is sharp in the sense
that it is not valid if $p=q$). We also investigate lineability of
the sets of norm-attaining and non-norm attaining operators.

\end{abstract}

\section{Introduction}

Let ${E}$ and ${F}$ be two Banach spaces. We denote by $\mathcal{L}(E,F)$ the
space of all bounded linear operators from $E$ into $F$. A linear operator
${T}\colon E\rightarrow F$ is said to attain its norm if there exists a $v\in
E$, $\Vert v\Vert_{E}=1$, such that
\[
\Vert{T}(v)\Vert_{F}=\Vert{T}\Vert_{\mathcal{L}(E,F)}:=\sup\limits_{e\in
\mathbb{S}_{E}}\Vert{T}(e)\Vert_{F},
\]
where $\mathbb{S}_{E}$ denotes the unit sphere of $E$, i.e., $\mathbb{S}%
_{E}=\{e\in E:\Vert e\Vert_{E}=1\}.$ We will denote by $\mathcal{N\!A}(E,F)$
the subset of all norm attaining bounded linear operators from $E$ into $F$.

The question whether a given linear operator attains its norm is doubtless one
of the most important lines of investigation from the applied functional
analysis point of view. Often, the solvability of certain (continuous or
discrete) differential equations is intrinsically related to the norm
attaining property of a determined linear operator acting between appropriate
Banach spaces.

When the target space is the real line, i.e., $F = \mathbb{R}$, a deep and, by
now, well known result due to James, see \cite{J}, asserts that
$\mathcal{N\!A}(E,\mathbb{R}) = \mathcal{L}(E,\mathbb{R})$ if and only if $E$
is reflexive. Another classical result in this theory, Bishop-Phelps' Theorem,
\cite{BPh}, states that $\mathcal{N\!A}(E,\mathbb{R})$ is always norm-dense in
$\mathcal{L}(E,\mathbb{R})$.

The question whether $\mathcal{N\!A}(E,F)$ is dense in $\mathcal{L}(E,F)$ for
an arbitrary Banach space $F$ becomes much more involved. A remarkable result
due to Lindenstrauss assures that if $F$ is reflexive then indeed
$\mathcal{N\!A}(E,F)$ is dense in $\mathcal{L}(E,F)$. This result was further
generalized by Bourgain in \cite{B}, who showed that the Radon-Nikodym
property on $F$ suffices for $\mathcal{N\!A}(E,F)$ to be dense in
$\mathcal{L}(E,F)$. On the converse, Gowers in \cite{G} showed there exists a
Banach space $\mathfrak{E}$, such that $\mathcal{N\!A}(\ell_{p},\mathfrak{E})$
is not dense in $\mathcal{L}(\ell_{p},\mathfrak{E})$, for $1<p<\infty$. The
case $p=1$ was settled by Acosta in \cite{Ac}.

The first goal of this note is to provide a simple yet useful characterization
of norm attaining operators acting on $\ell_{p}$ type spaces. Hereafter $T$
will always denote a bounded linear operator from $\ell_{p}$ into $\ell_{q}$.
For a not weakly null maximizing sequence for $T$ we mean a sequence $u^{n}%
\in\ell_{p}$, with, $\Vert u^{n}\Vert_{\ell_{p}}=1$, $\Vert T(u^{n}%
)\Vert_{\ell_{q}}\rightarrow\Vert T\Vert$ and $u^{n}$ does not converge weakly
to zero.

Initially, let us recall that Pitt's Theorem, \cite{P}, states that any
bounded linear operator $T\colon\ell_{p}\rightarrow\ell_{q}$, with $1\leq q<p$
is compact. Therefore, $\mathcal{N\!A}(\ell_{p};\ell_{q})=\mathcal{L}(\ell
_{p};\ell_{q})$ provided that $q<p$. On the other hand, for $1\leq p\leq
q<\infty$ it is well-known that $\mathcal{N\!A}(\ell_{p};\ell_{q}%
)\neq\mathcal{L}(\ell_{p};\ell_{q})$ (see \cite[Proposition 4.2]{J-Wolfe}).
Yet in the lights of Pitt's Theorem, for $p>q$, if $T\neq0$ and $(x^{n}%
)_{n=1}^{\infty}$ is a maximizing sequence for $T$, then $(x^{n}%
)_{n=1}^{\infty}$ is not weakly null. Clearly this result is no longer valid
for $p\leq q$ (the inclusion provides an example). Our first and main result
shows that when $p\leq q$, the existence of a not weakly null maximizing
sequence for $T$ occurs precisely when $T$ is norm attaining.

\begin{theorem}
\label{Att} Let $1<p<\infty$, $1\leq q<\infty$ and $T\colon\ell_{p}%
\rightarrow\ell_{q}$ be a bounded linear operator. Then $T$ attains its norm
if, and only if, there exists a not weakly null maximizing sequence for $T$.
\end{theorem}

Theorem \ref{Att} is sharp in the sense that it is no longer true for $p=1$.
In fact, if $1\leq p\leq q<\infty$ the operator $T\in\mathcal{L}(\ell_{p}%
;\ell_{q})$ given by
\begin{equation}
T(x)=\left(  \frac{nx_{n}}{n+1}\right)  _{n=1}^{\infty}\label{novo1}%
\end{equation}
does not attain its norm. The canonical basis $(e_{j})_{j=1}^{\infty}$ is a
maximizing sequence for $T$ which is not weakly null when $p=1$. The authors
thank R. Aron for this observation.

The strategy for proving Theorem \ref{Att} relies on an asymptotic analysis
involving weak convergence in $\ell_{p}$-type spaces. From the proof of
Theorem \ref{Att}, as long as $p\not =q$, we can actually infer
pre-compactness of any not weakly null maximizing sequence for a linear
operator $T\in\mathcal{L}(\ell_{p};\ell_{q})$. This is the content of our next result.

\begin{theorem}
\label{Comp. Max. Seq} Let $T\colon\ell_{p}\rightarrow\ell_{q}$ be a not
identically zero norm attaining operator and $1<p<\infty$, $1\leq q<\infty$,
$p\not =q$. Then any not weakly null maximizing sequence for $T$ has a norm
convergent subsequence.
\end{theorem}

Initially, let us point out that Theorem \ref{Comp. Max. Seq} is sharp in the
sense it does not hold true if $p=q$. Indeed, the sequence
\[
u^{n}:=\left(  \dfrac{1}{\sqrt[p]{2}},0,0,\cdots,\dfrac{1}{\sqrt[p]{2}%
},0,\cdots\right)  =\dfrac{1}{\sqrt[p]{2}}e_{1}+\dfrac{1}{\sqrt[p]{2}}e_{n}%
\]
is a not weakly null maximizing sequence for the identity map, $\text{Id}%
\colon\ell_{p}\rightarrow\ell_{p}$, but has no norm convergent subsequence.

Theorem \ref{Att} is particulary useful in discrete problems involving some
sort of symmetry or special invariances. For instance, in practical
applications, one is often able to find a hyperplane $\Pi:= \{ f(x)
\le\epsilon\}$ with $\epsilon> 0$ such that $\|T(\xi) \|_{\ell_{q}} < \|T\|$,
for all $\xi\in\Pi\cap B_{1}$. Thus a maximizing sequence can be found within
$B_{1} \cap\{ f(x) > \epsilon\}$. In particular such a maximizing sequence is
not weakly null.

The simplest norm-invariant operation for sequences is permutation. In the
sequel, we state a definition and afterwards a consequence of Theorem
\ref{Att} related to permutation of sequences.

\begin{definition}
Given a real sequence $\alpha=\{\alpha_{j}\}_{j=1}^{\infty}\in c_{0}$, the
non-increasing permutation of $\alpha$, denoted as $\sigma(\alpha)=\{\beta
_{j}\}_{j=1}^{\infty}$ is given by
\[
\beta_{1}:=\max\limits_{j\in\mathbb{N}}\{|\alpha_{j}|\},\quad\beta_{2}%
:=\max\limits_{j\in\mathbb{N}}\left(  \{|\alpha_{j}|\}\setminus\{\beta
_{1}\}\right)  ,\cdots\beta_{k}:=\max\limits_{j\in\mathbb{N}}\left(
\{|\alpha_{j}|\}\setminus\{\beta_{1},\beta_{2},\cdots,\beta_{k-1}\}\right)
,\cdots,
\]
where $\beta_{k}=0$ if $\left(  \{|\alpha_{j}|\}\setminus\{\beta_{1},\beta
_{2},\cdots,\beta_{k-1}\}\right)  =\phi$. A linear operator $T\colon\ell
_{p}\rightarrow\ell_{q}$ is said to be monotone with respect to non-increasing
permutation if $\Vert T(\sigma(x))\Vert_{\ell_{q}}\geq\Vert T(x)\Vert
_{\ell_{q}}$ for every $x\in\ell_{p}$.
\end{definition}

A typical, but not the only, way of verifying that a given operator $T$ is
monotone with respect to non-increasing permutation is by checking that
\[
\langle T e_{1}, e_{j} \rangle\ge\langle T e_{2}, e_{j} \rangle\ge\cdots
\ge\langle T e_{k}, e_{j} \rangle\ge\cdots\ge0, \quad\forall j=1,2,\cdots
\]
Concerning monotone with respect to non-increasing permutation operators, we
have the following general result.

\begin{theorem}
\label{monotone} Let $1<p<\infty$, $1\le q< \infty$ and $T\colon\ell_{p}
\to\ell_{q}$ be monotone with respect to non-increasing permutation. Assume
for some $\epsilon> 0$, $T\big
|_{\ell_{p+\epsilon}} \hookrightarrow\ell_{q} $ continuously. Then $T$ attains
its norm.
\end{theorem}


As another simple yet interesting application of Theorem \ref{Att} (and
Bessafa-Pelczy\'{n}ski selection principle) we obtain, up to subsequences, a
structural behavior of any maximizing sequence for an operator $T\in
\mathcal{L}(\ell_{p};\ell_{q})\setminus\mathcal{N\!A}(\ell_{p};\ell_{q})$. We
recall that if $\{e_{i}\}$ and $\{f_{i}\}$ are two basic sequences in Banach
spaces, then we say $\{e_{i}\}$ is equivalent to $\{f_{i}\}$ if for any
sequence of scalars $\{\lambda_{i}\}$, $\sum\limits_{i}\lambda_{i}e_{i}$
converges if and only if $\sum\limits_{i}\lambda_{i}f_{i}$ converges.

\begin{proposition}
\label{Coro} Let $1<p<\infty$, $1\leq q<\infty$ and $T\colon\ell
_{p}\rightarrow\ell_{q}$ be a bounded linear operator. Assume $T$ does not
attain its norm. Then, any maximizing sequence $u^{n}$ for $T$ has a
subsequence, $u^{n_{k}}$, such that $u^{n_{k}}$ is isometrically equivalent to
the canonical basis of $\ell_{p},$ whose image is equivalent to the canonical
basis of $\ell_{q}$.
\end{proposition}

In a parallel direction, we also carry out the investigation of lineability
properties related to the sets $\mathcal{N\!A}(X;\ell_{q})$ as well as
$\mathcal{L}(X;\ell_{q})\setminus\mathcal{N\!A}(X;\ell_{q})$. Recall that in
an infinite-dimensional vector space $X$, a set $A\subset X$ is said to be
lineable if $A\cup\{0\}$ contains an infinite-dimensional subspace. The term
\textquotedblleft lineable\textquotedblright\ seems to have been coined by V.
Gurariy and has been broadly explored in different contexts (see, for example
\cite{Aco, Ar, Ber} and references therein).

Let $X$ and $Y$ be Banach spaces. For a fixed vector $x_{0}\in\mathbb{S}_{X}$,
let us denote $\mathcal{N\!A}^{x_{0}}(X;Y)$ the set of all linear operators in
$\mathcal{L}(X;Y)$ that attain their norms at $x_{0}$, that is,
\[
\mathcal{N\!A}^{x_{0}}(X;Y):=\left\{  T\in\mathcal{N\!A}(X;Y):\Vert
Tx_{0}\Vert_{Y}=\Vert T\Vert_{\mathcal{L}(X;Y)}\right\}  .
\]
The linear structure of the sets $\mathcal{N\!A}(X;\mathbb{K})$ and
$\mathcal{L}(X;\mathbb{K})\smallsetminus\mathcal{N\!A}(X;\mathbb{K})$ was
subject of several recent works (see, e.g., \cite{Ac, Aco} and references
therein) and, of course, the geometry of $X$ plays a decisive role in this
study. If we replace $\mathbb{K}$ by an infinite-dimensional Banach space $Y$,
as it will be shown, it seems that the geometry of $Y$ will be decisive,
rather than the particular properties of $X$. We believe that the study is
lineability properties related to $\mathcal{N\!A}(X;Y)$, where $Y$ is a
hereditarily indecomposable space may be an interesting subject for further investigation.

It is worth mentioning that the presence of an infinite-dimensional Banach
space $Y$ in the place of the scalar field $\mathbb{K}$ allows to investigate
the lineability of sets of norm-attaining operators at a fixed point $x_{0}.$

In general, $\mathcal{N\!A}^{x_{0}}(X;Y)$ is a quite more restrictive subset
of $\mathcal{N\!A}(X;Y)$. Nevertheless we have managed to show that if $Y$
contains an isometric copy of $\ell_{q}$, then $\mathcal{N\!A}^{x_{0}}(X;Y)$
is lineable in $\mathcal{L}(X;Y)$. In particular $\mathcal{N\!A}^{x_{0}}%
(\ell_{p};\ell_{q})$ is lineable in $\mathcal{L}(\ell_{p};\ell_{q})$. This is
the content of our next result.

\begin{proposition}
\label{Lineability of NA} Let $X$ and $Y$ be Banach spaces so that $Y$
contains an isometric copy of $\ell_{q}$ for some $1\leq q<\infty,$ and let
$x_{0}\in\mathbb{S}_{X}.$ Then $\mathcal{N\!A}^{x_{0}}(X;Y)$ is lineable in
$\mathcal{L}(X;Y)$.\bigskip
\end{proposition}

An adaptation of the argument used to prove Proposition
\ref{Lineability of NA} allows us to conclude that $\mathcal{L}(\ell_{p}%
;\ell_{q})\smallsetminus\mathcal{N\!A}(\ell_{p};\ell_{q})$ is also lineable in
$\mathcal{L}(\ell_{p};\ell_{q})$ when $1\leq p\leq q<\infty$. In fact we prove
a more general result:

\begin{proposition}
\label{Lineability of NA22} Let $X$ and $Y$ be Banach spaces so that $Y$
contains an isometric copy of $\ell_{q}$ for some $1\leq q<\infty.$ If
$\mathcal{L}(X;\ell_{q})\neq\mathcal{N\!A}(X;\ell_{q}),$ then $\mathcal{L}%
(X;Y)\smallsetminus\mathcal{N\!A}(X;Y)$ is lineable in $\mathcal{L}(X;Y)$.
\end{proposition}

\begin{corollary}
\label{Lineability of L-NA} If $1\leq p\leq q<\infty$ then $\mathcal{L}%
(\ell_{p};\ell_{q})\smallsetminus\mathcal{N\!A}(\ell_{p};\ell_{q})$ is
lineable in $\mathcal{L}(\ell_{p};\ell_{q})$.
\end{corollary}

The arguments used throughout the article are fairly clear and simple in
nature. We do believe they provide insights for further generalizations to
more abstract settings.

\section{Characterization of operators in $\mathcal{N\!A}(\ell_{p};\ell_{q})$}

\textbf{Proof of Theorem \ref{Att}.} Clearly, if $T$ attains its norm, there
exists a not weakly null maximizing sequence for $T$. We shall just address
the converse. As mentioned in the introduction, when $p>q$, \textit{any}
bounded linear operator attains its norm. The really interesting situation for
us is therefore when $1<p\leq q$. From now on $e_{1},e_{2},...$ will denote
the canonical unit vectors of the sequence spaces.

Let $u^{n}$ be a not weakly null maximizing sequence for $T$. Passing to a
subsequence if necessary, we may assume,
\[
u^{n}\rightharpoonup u\not =0,
\]
where the symbol $\rightharpoonup$ stands for the weak convergence in
$\ell_{p}$. Our first observation is that there exists a subsequence $(v^{n})$
of $(u^{n})$ so that%
\begin{equation}
\left(  |\langle T(v^{n}),e_{i}\rangle-\langle T(u),e_{i}\rangle
|^{q-1}\right)  _{i=1}^{\infty}\rightharpoonup0 \label{nov21a}%
\end{equation}
in $\ell_{\frac{q}{q-1}}$. Indeed, since $\left(  |\langle T(u^{n}%
),e_{i}\rangle-\langle T(u),e_{i}\rangle|^{q-1}\right)  _{i=1}^{\infty}$ is a
bounded sequence in $\ell_{\frac{q}{q-1}}$, there is a subsequence $(v^{n})$
of $(u^{n})$ so that%
\[
\left(  |\langle T(v^{n}),e_{i}\rangle-\langle T(u),e_{i}\rangle
|^{q-1}\right)  _{i=1}^{\infty}\rightharpoonup f,
\]
for some $f\in\ell_{\frac{q}{q-1}}$. Since $v^{n}\rightharpoonup u$ we have
that $T(v^{n})\rightharpoonup T(u);$ so, for each $i=1,2,\cdots$
\[
\langle T(v^{n}),e_{i}\rangle-\langle T(u),e_{i}\rangle\rightarrow0,
\]
thus $f$ must be the null vector.\bigskip

Let $r>1$ be a real number and let us consider the auxiliary function
$\varphi_{r}:\mathbb{R}\setminus\{1\}\rightarrow\mathbb{R}$ given by
\[
\varphi_{r}(X)=\dfrac{\Big |~|X|^{r}-|X-1|^{r}-1~\Big |}{|X-1|^{r-1}}.
\]
It is simple to verify that
\[
\lim\limits_{|X|\rightarrow\infty}\varphi(X)=r.
\]
Hence, given $\varepsilon>0$, we can find a constant $C_{\varepsilon}$ such
that
\[
\Big |~|X|^{r}-|X-1|^{r}-1~\Big |\leq C_{\varepsilon}|X-1|^{r-1}%
\]
whenever $|X-1|>\varepsilon$. On the other hand if $|X-1|\leq\varepsilon$, we
have
\[
\Big |~|X|^{r}-|X-1|^{r}-1~\Big |\leq\varepsilon^{r}+\tilde{\delta
}(\varepsilon),
\]
where
\[
\tilde{\delta}(\varepsilon):=\sup\limits_{|t-1|\leq\varepsilon}|~|t|^{r}-1~|.
\]
Adding up the above two inequalities we obtain
\begin{equation}
\Big |~|X|^{r}-|X-1|^{r}-1~\Big |\leq C_{\varepsilon}|X-1|^{r-1}%
+\delta(\varepsilon), \label{RealAn-Inq}%
\end{equation}
for every $X\in\mathbb{R}$, where $\delta$ is a modulus of continuity.

The idea now is to apply estimate (\ref{RealAn-Inq}) to each coordinate of
$T(v^{n})$ in $\ell_{q}$. More precisely, we apply inequality
(\ref{RealAn-Inq}) to $r=q$ and
\[
X_{i}:=\dfrac{\langle T(v^{n}),e_{i}\rangle}{\langle T(u),e_{i}\rangle},
\]
whenever $\langle T(u),e_{i}\rangle$ is nonzero. In any case, when we sum
these inequalities up, we obtain
\begin{equation}%
\begin{array}
[c]{lll}%
\sum\limits_{i=1}^{\infty}\Big ||\langle T(v^{n}),e_{i}\rangle|^{q}-|\langle
T(v^{n}),e_{i}\rangle-\langle T(u),e_{i}\rangle|^{q}-|\langle T(u),e_{i}%
\rangle|^{q}\Big | & \leq & C_{\varepsilon}\Delta_{n}\\
& + & \delta(\varepsilon)\Vert T(u)\Vert_{\ell_{q}}^{q},
\end{array}
\label{BoundBL}%
\end{equation}
where
\[
\Delta_{n}=\sum\limits_{i=1}^{\infty}|\langle T(u),e_{i}\rangle|\cdot|\langle
T(v^{n}),e_{i}\rangle-\langle T(u),e_{i}\rangle|^{q-1}.
\]
Since $T(u)\in\ell_{q}=\big [\ell_{\frac{q}{q-1}}\big ]^{\ast}$, we have
\[
\Delta_{n}=\langle T(u),\left(  |\langle T(v^{n}),e_{i}\rangle-\langle
T(u),e_{i}\rangle|^{q-1}\right)  _{i=1}^{\infty}\rangle
\]
and using (\ref{nov21a}) it follows that
\begin{equation}
\lim\limits_{n\rightarrow\infty}\Delta_{n}=\langle T(u),0\rangle=0.
\label{error}%
\end{equation}
Letting $n\rightarrow\infty$ in (\ref{BoundBL}) we get, for every
$\varepsilon>0$,
\begin{align*}
&  \limsup_{n\rightarrow\infty}\Big |\sum\limits_{i=1}^{\infty}\left(
|\langle T(v^{n}),e_{i}\rangle|^{q}-|\langle T(v^{n}),e_{i}\rangle-\langle
T(u),e_{i}\rangle|^{q}-|\langle T(u),e_{i}\rangle|^{q}\right)  \Big |\\
&  \leq\delta(\varepsilon)\Vert T(u)\Vert_{q}^{q},
\end{align*}

Making $\varepsilon\searrow0$ we conclude that the $\limsup$ is in fact the
limit and is equal to zero:%
\[
\lim_{n\rightarrow\infty}\Big |\sum\limits_{i=1}^{\infty}\left(  |\langle
T(v^{n}),e_{i}\rangle|^{q}-|\langle T(v^{n}),e_{i}\rangle-\langle
T(u),e_{i}\rangle|^{q}-|\langle T(u),e_{i}\rangle|^{q}\right)  \Big |=0.
\]
In particular,
\begin{equation}
\Vert T(v^{n})\Vert_{\ell_{q}}^{q}=\Big (\Vert Tu\Vert_{\ell_{q}}^{q}+\Vert
T(v^{n}-u)\Vert_{\ell_{q}}^{q}+\text{o}(1)\Big ), \label{Teix-L}%
\end{equation}
where $\text{o}(1)\rightarrow0$ as $n\rightarrow\infty$.

A similar computation, using $r=p$ on inequality (\ref{RealAn-Inq}) and
applying it on the points
\[
Y_{i}:=\dfrac{\langle v^{n},e_{i}\rangle}{\langle u,e_{i}\rangle}%
\]
can be performed, as long as $\langle u,e_{i}\rangle$ is nonzero. Reasoning as
before, we reach
\begin{equation}
\Vert w^{n}-u\Vert_{\ell_{p}}^{p}=1-\Vert u\Vert_{\ell_{p}}^{p}+\text{o}(1)
\label{Teix-L2}%
\end{equation}
for some subsequence $(w^{n})$ of $(v^{n}).$ Combining (\ref{Teix-L}) and
(\ref{Teix-L2}) with the well known inequality
\[
(\alpha+\beta)^{\theta}\leq\alpha^{\theta}+\beta^{\theta},
\]
for $\alpha,\beta\geq0$ and $0\leq\theta\leq1$, we can write
\begin{equation}%
\begin{array}
[c]{lll}%
\Vert T(w^{n})\Vert_{\ell_{q}}^{p} & \leq & \Big (\Vert T(u)\Vert_{\ell_{q}%
}^{q}+\Vert T(w^{n}-u)\Vert_{\ell_{q}}^{q}+\text{o}(1)\Big )^{p/q}\\
& \leq & \Vert T(u)\Vert_{\ell_{q}}^{p}+\Vert T(w^{n}-u)\Vert_{\ell_{q}}%
^{p}+\text{o}(1)\\
& \leq & \Vert T(u)\Vert_{\ell_{q}}^{p}+\Vert T\Vert^{p}\cdot\Vert
w^{n}-u\Vert_{\ell_{p}}^{p}+\text{o}(1)\\
& \leq & \Vert T(u)\Vert_{\ell_{q}}^{p}+\Vert T\Vert^{p}\cdot\lbrack1-\Vert
u\Vert_{\ell_{p}}^{p}+\text{o}(1)]+\text{o}(1).
\end{array}
\label{jan28-09}%
\end{equation}
Letting $n\rightarrow\infty$, we finally obtain
\[
\Vert T(u)\Vert_{\ell_{q}}\geq\Vert T\Vert\cdot\Vert u\Vert_{\ell_{p}},
\]
which finishes the proof of Theorem \ref{Att}. \hfill$\square$

\section{A disguise of Theorem \ref{Att}: Pre-compactness of maximizing
sequences}


\textbf{Proof of Theorem \ref{Comp. Max. Seq}.} Let $T\colon\ell
_{p}\rightarrow\ell_{q}$, with $T\not =0$ be a norm attaining operator and
$x^{n}$ a maximizing sequence on $\mathbb{S}_{\ell_{p}}$ that does not
converge weakly to zero. Up to a subsequence, $x^{n}$ converges weakly to a
point $x_{0}$. By uniform convexity of $\ell_{p}$ (or if you prefer, equation
(\ref{Teix-L2})) it suffices to show $x_{0}\in\mathbb{S}_{\ell_{p}}$.

When $p>q$, our thesis is a consequence of Pitt's Theorem. Indeed, since $T$
is a compact operator, $T(x^{n})$ converges strongly to $T(x_{0})$ in
$\ell_{q}$.

Since $x^{n}\rightharpoonup x_{0}$ it follows that $T(x^{n})\rightarrow
T(x_{0})$ and $\left\Vert x_{0}\right\Vert \leq1.$ But $\Vert T(x^{n}%
)\Vert_{\ell_{q}}\rightarrow\Vert T\Vert;$ so we conclude that
\[
\Vert T\Vert=\Vert T(x_{0})\Vert_{\ell_{q}}%
\]
and hence $\Vert x_{0}\Vert_{\ell_{p}}=1$.

When $p<q$, we argue as follows: we may assume $x_{0}\not =0$. As before, we
have to verify that $x_{0}\in\mathbb{S}_{\ell_{p}}$. Reasoning as in the proof
of Theorem \ref{Att} (see (\ref{jan28-09})), we can write
\begin{equation}
\Vert T(x^{n})\Vert_{\ell_{q}}^{q}\leq\Vert T(x_{0})\Vert_{\ell_{q}}^{q}+\Vert
T\Vert^{q}\left(  1-\Vert x_{0}\Vert_{\ell_{p}}^{p}\right)  ^{q/p}%
+\text{o}(1). \label{pre-comp. Eq01}%
\end{equation}
Since $\Vert T(x^{n})\Vert_{\ell_{q}}=\Vert T\Vert+\text{o}(1)$ and (from the
proof of Theorem \textbf{\ref{Att})} $\Vert T(x_{0})\Vert_{\ell_{q}}=\Vert
T\Vert\cdot\Vert x_{0}\Vert_{\ell_{p}}$, equation (\ref{pre-comp. Eq01}) leads
to
\begin{equation}
1\leq\Vert x_{0}\Vert_{\ell_{p}}^{q}+\left(  1-\Vert x_{0}\Vert_{\ell_{p}}%
^{p}\right)  ^{q/p}. \label{pre-comp. Eq02}%
\end{equation}
Finally, since $q/p>1$, equation (\ref{pre-comp. Eq02}) implies $1-\Vert
x_{0}\Vert^{p}=0$; otherwise the strict inequality would hold
\[
1=\left[  \Vert x_{0}\Vert_{\ell_{p}}^{p}+(1-\Vert x_{0}\Vert_{\ell_{p}}%
^{p})\right]  ^{q/p}>\Vert x_{0}\Vert_{\ell_{p}}^{q}+\left(  1-\Vert
x_{0}\Vert_{\ell_{p}}^{p}\right)  ^{q/p}.
\]
which drives us to a contradiction. \hfill$\square$


\section{Two Applications of Theorem \ref{Att}}

\textbf{Proof of Theorem \ref{monotone}.} Let $T\colon\ell_{p}\rightarrow
\ell_{q}$ be monotone with respect to non-increasing permutation and consider
$x^{n}$ a maximizing sequence for $T$. We may and will assume $T$ is not the
zero operator. Since $\Vert T(\sigma(x^{n}))\Vert_{\ell_{q}}\geq\Vert
T(x^{n})\Vert_{\ell_{q}}$, and $\Vert\sigma(x^{n})\Vert_{\ell_{p}}=\Vert
x^{n}\Vert_{\ell_{p}}=1$, $y^{n}:=\sigma(x^{n})$ is too a maximizing sequence
for $T$. In view of Theorem \ref{Att} it suffices to verify $y^{n}$ is not
weakly null. For that we argue as follows: suppose, for sake of contradiction,
that $y^{n}$ does converge weakly to zero. Since $y^{n}$ is in non-increasing
order, it would imply $\Vert y^{n}\Vert_{\ell_{\infty}}=\text{o}(1)$ as
$n\rightarrow\infty$, and therefore
\begin{equation}
\Vert y^{n}\Vert_{\ell_{p+\epsilon}}\leq\Vert y^{n}\Vert_{\ell_{\infty}%
}^{\frac{\epsilon}{p+\epsilon}}\cdot\Vert y^{n}\Vert_{\ell_{p}}^{\frac
{p}{p+\epsilon}}=\text{o}(1)\text{ as }n\rightarrow\infty.
\label{monotone Eq0}%
\end{equation}
By continuity, (\ref{monotone Eq0}) would lead us to
\[
\Vert T\Vert=\lim\limits_{n\rightarrow\infty}\Vert Ty^{n}\Vert_{q}=0,
\]
which is a contradiction to our earlier assumption, $T\not \equiv 0$.
\hfill$\square$ \bigskip

\noindent\textbf{Proof of Proposition \ref{Coro}.} Assume $T\colon\ell
_{p}\rightarrow\ell_{q}$ does not attain its norm. From Theorem \ref{Att}, for
any maximizing sequence $u^{n}$, one has
\[
u^{n}\rightharpoonup0\text{ in }\ell_{p}.
\]
Therefore, because of Bessaga-Pelczy\'{n}ski selection principle, see, e.g.,
\cite{LiT}, there exists a infinite subset of the natural numbers,
$\mathbb{N}_{1}\subseteq\mathbb{N}$, such that $\left\{  u^{n}\right\}
_{n\in\mathbb{N}_{1}}$ is a basic sequence equivalent to a block basic
sequence of the canonical basis of $\ell_{p}$. But now it is simple to show
that actually
\[
\left\{  u^{n}\right\}  _{n\in\mathbb{N}_{1}}\text{ is equivalent to the
canonical basis of }\ell_{p}.
\]
Furthermore, $\overline{\text{span}}\left\{  u^{n}\right\}  _{n\in
\mathbb{N}_{1}}$ is isometric to $\ell_{p}$. Indeed, because $\Vert u^{n}%
\Vert_{p}=1$ and $\left\{  u^{n}\right\}  _{n\in\mathbb{N}_{1}}$ is equivalent
to a block basic sequence of $\{e_{i}\}$, we can find scalars $\gamma_{i}$
such that
\[
u^{i}=\sum\limits_{k=r_{i}+1}^{r_{i+1}}\gamma_{k}e_{k},~\forall i\in
\mathbb{N}_{1},
\]
with
\[
\sum\limits_{k=r_{i}+1}^{r_{i+1}}|\gamma_{k}|^{p}=1,~\forall i\in
\mathbb{N}_{1}.
\]
Now,
\[%
\begin{array}
[c]{lll}%
\left\Vert \sum\limits_{i=1}^{M}a_{i}u^{i}\right\Vert _{p} & = & \left(
\sum\limits_{i=1}^{M}\sum\limits_{k=r_{i}+1}^{r_{i+1}}|a_{i}|^{p}|\gamma
_{k}|^{p}\right)  ^{1/p}\\
& = & \left(  \sum\limits_{i=1}^{M}|a_{i}|^{p}\sum\limits_{k=r_{i}+1}%
^{r_{i+1}}|\gamma_{k}|^{p}\right)  ^{1/p}\\
& = & \left\Vert \sum\limits_{i=1}^{M}a_{i}e_{i}\right\Vert _{p}.
\end{array}
\]
Now, as long as $T\not \equiv 0$, the sequence $\left\{  T(u^{n})\right\}
_{n\in\mathbb{N}_{1}}$ is weakly null but
\[
\Vert T(u^{n})\Vert_{q}\rightarrow\Vert T\Vert\not =0.
\]
Thus, applying the same argument as before to the sequence $\left\{
T(u^{n})\right\}  _{n\in\mathbb{N}_{1}}$, we find a $\mathbb{N}_{2}%
\subseteq\mathbb{N}_{1}$, such that the sequence
\[
\left\{  T(u^{n})\right\}  _{n\in\mathbb{N}_{2}}\text{is equivalent to the
canonical basis of }\ell_{q}%
\]
and the proof of the Corollary is complete. \hfill$\square$

\section{Lineability of the set of norm attaining operators at a fixed point}


\textbf{Proof of Proposition \ref{Lineability of NA}.} Our first observation
is that it suffices to prove Proposition \ref{Lineability of NA} for
$Y=\ell_{q}$. Using Hahn-Banach Theorem it is not difficult to show that
$\mathcal{N\!A}^{x_{0}}(X;\ell_{q})\neq\{0\}.$

Hereafter we fix a nonzero element $u\in\mathcal{N\!A}^{x_{0}}(X;\ell_{q}).$
We can write the set of positive integers $\mathbb{N}$ as
\[
\mathbb{N=}%
{\displaystyle\bigcup\limits_{k=1}^{\infty}}
A_{k},
\]
where each
\begin{equation}
A_{k}:=\{a_{1}^{(k)}<a_{2}^{(k)}<...\}\label{hhggpp}%
\end{equation}
has the same cardinality as $\mathbb{N}$ and the sets $A_{k}$ are pairwise
disjoint. For each positive integer $k$, we define%
\[
\ell_{q}^{(k)}:=\left\{  x\in\ell_{q}:x_{j}=0\text{ if }j\notin A_{k}\right\}
.
\]
In the sequel, for each $k$ fixed, let $u_{k}:X\rightarrow\ell_{q}^{(k)}$ be
given by
\[
\left(  u_{k}(x)\right)  _{a_{j}^{(k)}}=\left(  u(x)\right)  _{j},\quad\forall
j,k\in\mathbb{N}.
\]
Finally, for $k$ fixed, let $v_{k}:X\rightarrow\ell_{q}$ be given by%
\[
v_{k}=i_{k}\circ u_{k},
\]
where $i_{k}:\ell_{q}^{(k)}\rightarrow\ell_{q}$ is the canonical inclusion.
Note that
\[
\left\Vert v_{k}(x)\right\Vert =\left\Vert u_{k}(x)\right\Vert =\left\Vert
u(x)\right\Vert
\]
for every positive integer $k$ and $x\in X.$ Thus, each $v_{k}$ attains its
norm at $x_{0}$. From the fact that the operators $v_{k}$ have disjoint
supports it is easy to verify that
\[
\{v_{1},v_{2},...\}
\]
is a linearly independent set. It just remains to verify that any operator in
\[
\text{span}\{v_{1},v_{2},...\}
\]
attains its norm at $x_{0}.$ For notation convenience, we will show that
$av_{1}+bv_{2}$ attains its norm at $x_{0}$ for any choice of scalars $a,b$.
We compute%

\begin{align*}
\left\Vert av_{1}+bv_{2}\right\Vert ^{q}  &  =\sup_{\left\Vert x\right\Vert
\leq1}\left\Vert av_{1}(x)+bv_{2}(x)\right\Vert ^{q}\\
&  \overset{(\ast)}{=}\sup_{\left\Vert x\right\Vert \leq1}%
{\displaystyle\sum\limits_{k}}
\left\vert a(v_{1}(x))_{k}\right\vert ^{q}+%
{\displaystyle\sum\limits_{k}}
\left\vert b(v_{2}(x))_{k}\right\vert ^{q}\\
&  =\sup_{\left\Vert x\right\Vert \leq1}\left\vert a\right\vert ^{q}%
{\displaystyle\sum\limits_{k}}
\left\vert (v_{1}(x))_{k}\right\vert ^{q}+\left\vert b\right\vert ^{q}%
{\displaystyle\sum\limits_{k}}
\left\vert (v_{2}(x))_{k}\right\vert ^{q}\\
&  =\left\vert a\right\vert ^{q}%
{\displaystyle\sum\limits_{k}}
\left\vert (v_{1}(x_{0}))_{k}\right\vert ^{q}+\left\vert b\right\vert ^{q}%
{\displaystyle\sum\limits_{k}}
\left\vert (v_{2}(x_{0}))_{k}\right\vert ^{q}\\
&  =\left\Vert av_{1}(x_{0})+bv_{2}(x_{0})\right\Vert ^{q}.
\end{align*}
Thus, indeed $av_{1} + bv_{2}$ attains its norm at $x_{0}$. Equality $(*)$
holds since $v_{1}$ and $v_{2}$ have disjoint supports. \hfill$\square$

\section{Lineability of sets of non-norm-attaining operators}

\textbf{Proof of Proposition \ref{Lineability of NA22}.} We just need to deal
with the case $Y=\ell_{q}$. Let $T:X\rightarrow\ell_{q}$ be a
non-norm-attaining operator. Again, we write $\mathbb{N}$ as
\[
\mathbb{N=}%
{\displaystyle\bigcup\limits_{k=1}^{\infty}}
A_{k},
\]
with the $A_{k}$ as in (\ref{hhggpp}).Again, for each positive integer $k$,
let%
\[
\ell_{q}^{(k)}:=\left\{  x\in\ell_{q};x_{j}=0\text{ if }j\notin A_{k}\right\}
.
\]
For each $k$, we consider $T_{k}:X\rightarrow\ell_{q}^{(k)}$ defined as
\[
\left(  T_{k}(x)\right)  _{a_{j}^{(k)}}=\left(  T(x)\right)  _{j},\quad\forall
j,k\in\mathbb{N}.
\]
For every $k$, let $v_{k}:X\rightarrow\ell_{q}$ given by%
\[
v_{k}=i_{k}\circ T_{k},
\]
where $i_{k}:\ell_{q}^{(k)}\rightarrow\ell_{q}$ is the canonical inclusion.
So, as in the previous proof, each $v_{k}$ does not attain its norm and
\[
\{v_{1},v_{2},...\}
\]
is a linearly independent set. It remains to be shown that any operator in
\[
\text{span}\{v_{1},v_{2},...\}
\]
does not attain their norms. Again, for notation convenience, let us restrict
our computation to $av_{1}+bv_{2}$, for any choice of scalars $a,b$ (of
course, at least one of them is chosen to be different from zero). To show
that $av_{1}+bv_{2}$ does not attain its norm we argue as follows:
\begin{align*}
\left\Vert av_{1}+bv_{2}\right\Vert ^{q} &  =\sup_{\left\Vert x\right\Vert
\leq1}\left\Vert av_{1}(x)+bv_{2}(x)\right\Vert ^{q}\\
&  =\sup_{\left\Vert x\right\Vert \leq1}%
{\displaystyle\sum\limits_{k}}
\left(  \left\vert a(v_{1}(x))_{k}\right\vert ^{q}+%
{\displaystyle\sum\limits_{k}}
\left\vert b(v_{2}(x))_{k}\right\vert ^{q}\right)  \\
&  \leq\left\vert a\right\vert ^{q}\left\Vert v_{1}\right\Vert ^{q}+\left\vert
b\right\vert ^{q}\left\Vert v_{2}\right\Vert ^{q}.
\end{align*}
On the other hand, for every natural number $n$ and $\varepsilon=\frac{1}{n}$
we can find $x_{n}\in\mathbb{S}_{X}$ so that%
\[
\left\Vert v_{j}(x_{n})\right\Vert \geq\left\Vert v_{j}\right\Vert
-\varepsilon,\text{ }j=1,2
\]
and hence%
\begin{align*}
\left\Vert (av_{1}+bv_{2})(x_{n})\right\Vert ^{q} &  =\left\vert a\right\vert
^{q}\left\Vert v_{1}(x_{n})\right\Vert ^{q}+\left\vert b\right\vert
^{q}\left\Vert v_{2}(x_{n})\right\Vert ^{q}\\
&  \geq\left\vert a\right\vert ^{q}(\left\Vert v_{1}\right\Vert -\varepsilon
)^{q}+\left\vert b\right\vert ^{q}(\left\Vert v_{2}\right\Vert -\varepsilon
)^{q}.
\end{align*}
So we conclude that%
\[
\left\Vert av_{1}+bv_{2}\right\Vert ^{q}=\left\vert a\right\vert
^{q}\left\Vert v_{1}\right\Vert ^{q}+\left\vert b\right\vert ^{q}\left\Vert
v_{2}\right\Vert ^{q}.
\]

Besides, if $\left\Vert x\right\Vert _{X}=1$, since $v_{1}$ and $v_{2}$ do not
attain their norms, we get%
\begin{align*}
\left\Vert (av_{1}+bv_{2})(x)\right\Vert ^{q}  &  =\left\vert a\right\vert
^{q}\left\Vert v_{1}(x)\right\Vert ^{q}+\left\vert b\right\vert ^{q}\left\Vert
v_{2}(x)\right\Vert ^{q}\\
&  <\left\vert a\right\vert ^{q}\left\Vert v_{1}\right\Vert ^{q}+\left\vert
b\right\vert ^{q}\left\Vert v_{2}\right\Vert ^{q}\\
&  =\left\Vert av_{1}+bv_{2}\right\Vert ^{q}.
\end{align*}
We conclude that $av_{1}+bv_{2}$ belongs to $\mathcal{L}(X;\ell_{q}%
)\smallsetminus\mathcal{N\!A}(X;\ell_{q})$. The general case is similar. The
proof of Proposition \ref{Lineability of NA22} is completed. \hfill$\square$

\vfill

\noindent\textsc{Daniel Pellegrino} \hfill\textsc{Eduardo V. Teixeira}\newline
Universidade Federal da Para\'{\i}ba \hfill Univerisdade Federal do
Cear\'{a}\newline Departamento de Matem\'{a}tica\hfill Departamento de
Matem\'{a}tica\newline Jo\~{a}o Pessoa-PB, Brazil, CEP 58.051-900\hfill Av.
Humberto Monte, s/n\newline E-mail: \texttt{dmpellegrino@gmail.com} \hfill
Fortaleza-CE, Brazil. CEP 60.455-760\newline.\hfill E-mail:
\texttt{eteixeira@pq.cnpq.br}

\end{document}